\documentstyle{amsart}
\theoremstyle{theorem}
\newtheorem{defn}{Definition}[section]
\newtheorem{thm}[defn]{Theorem}
\newtheorem{ex}[defn]{Example}
\newtheorem{lem}[defn]{Lemma}

\newtheorem{cor}[defn]{Corollary}

\newtheorem{prob}{Problem}

\def\R{{\mathbb R}}

\def\N{{\mathbb  N}}

\def\e{\epsilon}

\def\ss{\subseteq}

\def\hm{{\mathcal HM}}
\title{On Haar Meager Sets}
\author[U. B. Darji]{Udayan B. Darji}
\address{Department of Mathematics University of Louisville
Louisville, KY 40292,USA} 
\email{ubdarj01@@louisville.edu}
\begin{document}
\begin{abstract} The notion of Haar null set was introduced by J. P. R. Christensen in 1973 and reintroduced in 1992 in the context of dynamical 
systems by Hunt, Sauer and Yorke.  During the last twenty years this notion has been useful in studying exceptional sets in diverse areas.
These include  analysis, dynamical systems, group theory, and descriptive set theory. Inspired by these various results, we introduce the topological analogue of 
the notion of Haar
null set. We call it Haar meager set. We prove some basic properties of this notion, state some open problems and suggest a possible line of 
investigation which may lead to the unification of these two notions in certain context.
\end{abstract} 
\maketitle 
\section{introduction}
Often in various branches of mathematics one would like to show that a certain phenomenon occurs very rarely or occurs rather 
frequently. This is particularly the case in analysis and probability theory. Hence we have the notion of measure zero set and
the notion of almost everywhere or  almost sure. Each locally compact group admits a translation invariant regular Borel measure
which is finite on compact sets and positive on open sets. Such a measure is unique, up to multiplication by a constant, and is 
called Haar measure. More precisely, if $\mu$ and $\nu$ are two such measure and $\nu$ is not identically zero,
then $\mu = c \nu$ for some $c \ge 0$. Of course, in $\R ^n$ this is the Lebesgue measure. If one has no prior bias towards any particular set of 
points one uses this measure. However, Haar measures only exist on locally compact groups. Therefore, in many important
spaces such as $C([0,1])$, the space of continuos real-valued functions on $[0,1]$ and $S_{\infty}$, the group of permutations
on $\N$, one does not have a suitable notions of smallness. As these objects are complete topological groups with natural topologies,
one can use the notion of meagerness as the notion of smallness. One of the earliest instance of such a result is a 
theorem of Banach which says that the set of all functions in $C([0,1])$ which are differentiable at some point is meager in
$C([0,1])$. 

The notion of meagerness is a topological one. Often it fails to capture the essence of certain properties studied in analysis, 
dynamical systems, etc. To remedy this, J. P. R. Christensen in 1973 \cite{christensen} introduced what he called ``Haar null" sets. The beauty of this concept
is that in locally compact group it is equivalent to the notion of measure zero set under the Haar measure and at the same time it is meaningful
in Polish groups in general.  Unaware of the result of Christensen, Hunt, Sauer and Yorke \cite{hsy} found this notion essential in 
the context of dynamical systems, reintroduced it and gave some applications of it in dynamics. Since then this notion has 
found applications in many diverse places. We only list a handful here. For applications in dynamics see \cite{hsy}, \cite{oy}. Its relations
to descriptive set theory and group theory can be found in \cite{dm}, \cite{solecki}, \cite{do}. Some applications of this notion
in real analysis can be found in \cite{ds},\cite{hunt}, \cite{kolar}, \cite{z}. 

Throughout $X$ is a fixed abelian Polish group. Following Christensen \cite{christensen} we say that $A \subseteq X$ is {\bf Haar null} if there is a Borel probability measure  $\mu$ on $X$ and a Borel set $B \subseteq X$
such that $A \subseteq B$ and $\mu(x +B )=0$ for all $x \in X$. By using standard tools of Borel measures on Polish groups, it can be shown that $\mu$ can be chosen so that its support is a Cantor set, i.e., a $0$-dimensional, compact metric space with no isolated points. In the original paper, Christensen showed that the set of Haar null sets forms a $\sigma$-ideal and in locally compact groups the notion of Haar null is equivalent to having measure zero under the Haar measure.
Often in the literature Haar null sets are called {\bf shy sets} and the complement of Haar null sets are called {\bf prevalent sets}.

There are explicit examples of Haar null sets which are also meager sets. For example, Hunt \cite{hunt} showed that the set of  functions in $C([0,1])$ which are differentiable at some point is Haar null, complementing the classical result of Banach that this set is meager. Kolar \cite{kolar}
introduced a notion of small sets  which he called {\bf HP-small}. This notion captures the concept of Haar null as well as  $\sigma$-porosity. (All $\sigma$-porous sets are meager. See \cite{z2} for more on $\sigma$-porosity.) Recently, Elekes and Stepr\=ans investigated some deep foundational properties of Haar null sets and analogous problems concerning meager sets. We refer the interested reader to their paper \cite{es}.

In this note we introduce what we call Haar meager sets. The purpose of this definition is to have a topological notion of smallness which is analogous to Christensen's measure theoretic notion of smallness. We show that Haar meager sets form a $\sigma$-ideal. We show every Haar meager set is meager. Example \ref{ex1} shows that some type of definability is necessary in the definition of Haar null to obtain this result.  We next show that in locally compact groups the notions  Haar meagerness and meagerness coincide. Using a result of Solecki, we observe that every non locally compact Polish group admits a closed nowhere dense set which is not Haar meager. We also formulate a criteria equivalent to Haar meagerness which may be useful in applications. 

In the next section we prove the main results. We end this note in the third section by stating some problems and suggesting some line of investigation. 

\section{main results}
Throughout, $X$ is an abelian Polish group.
\begin{defn} A set $A \subseteq X$ is {\bf Haar meager}
means that
there is a Borel set $B$ with $A \subseteq B$, a compact metric space $K$ and a continuous function $f: K \rightarrow X$
such that $f^{-1}(B +x)$ is meager
in $K$ for all $x \in X$. We call $f$
{\em a witness function} for $A$. The collection of all Haar meager sets is denoted by $\hm$.  \end{defn}
\begin{thm}\label{haarimpmeager} Let $A \ss X$ be Haar meager. Then, $A$ is meager. \end{thm}
\begin{pf} Since subsets of meager sets are meager, it suffices to prove the theorem for sets
which are Borel and Haar meager. Let $A$ be a Borel Haar meager set.  Let $K$ be a compact set and $f: K \rightarrow X$ be a
continuous function which witnesses that $A$ is Haar meager.
Let $\Sigma _A =\{(x,y) \in X \times K: f(y) \in (-x+ A) \}$. As $A$ is Borel, $\Sigma _A$ is Borel
and hence has the Baire
property. Also, for each $x \in X$, $\{y \in K: (x,y) \in \Sigma_ A \}$ is meager
in $K$ since $A$ is Haar meager. By Kuratowski-Ulam Theorem $\Sigma_A$ is meager in $X \times K$.
Again by Kuratowski-Ulam theorem we may choose $y \in K$ such that $B=\{x: (x,y) \in \Sigma_A\} $ is
meager in $X$. Note that $B=\{x: f(y) \in (-x+A)\} =\{x: x \in (-f(y)+A)\} = -f(y)+A$. Hence $-f(y)+A$
is meager in $X$. Since meager sets are translation invariant in a Polish group, we have that $A$ is
meager in $X$. \end{pf}
The following example shows that some sort of definability condition is necessary in
the definition of Haar meager sets in order for Haar meager sets to meager.
\begin{ex}\label{ex1} (Assume the Continuum Hypothesis) There is a subset $A \ss \mathbb R$ and a compact set $K \ss \mathbb R$ such
that $(A +x) \cap K $ is countable for all $x \in \mathbb R$, yet $A$ is not meager.  \end{ex}
\begin{pf} Let $K$ be the standard middle third Cantor set. Let $\{K_{\alpha}\}_{\alpha < \omega _1}$
be an enumeration of all translates of $K$ and let $\{C_{\alpha}\}_{\alpha < \omega _1}$ be an enumeration
of closed nowhere dense subsets of $\mathbb R$. For each $\alpha < \omega _1$, let $p_{\alpha} \in {\mathbb R}
 \setminus \bigcup_{\beta \le \alpha} \left ( K_{\beta} \cup C_{\beta} \right ) $. Let $A=\{p_{\alpha}: \alpha < \omega _1\}$.
 Then, $A$ is the desired set. \end{pf}
\begin{thm} If $X$ is locally compact and $A \subset X$ is meager, then $A$ is Haar meager. \end{thm}
\begin{pf}  Let $U \ss X$ be an
open set such that $\overline {U}$ is compact. That
$A$ is Haar meager is witnessed by $\overline{U}$ and the identity function on $\overline{U}$. \end{pf}
\begin{cor}\label{meagereqhaar} In a locally compact group $X$ a set is Haar meager iff it is meager. \end{cor}
\begin{ex} Every nonlocally compact Polish group $X$ has a closed meager set which is not Haar meager. \end{ex}
\begin{pf} This simply follows from the following result of Solecki.
\begin{thm} (Solecki \cite{solecki}) Let $X$ be a nonlocally compact Polish group. Then, there is a closed set $F \ss X$ and a continuous function
$f:F \rightarrow \{0,1\}^{\omega}$ such that each $f^{-1}(x)$, $x\in \{0,1\}^{\omega}$, contains a translate of every compact
subset of $X$. \end{thm}
Since $X$ is separable, there is $x \in \{0,1\} ^{\omega}$ such that $f^{-1}(x)$ is nowhere dense. This set is clearly not Haar
meager. \end{pf}
Our next goal is to prove that collection of Haar meager sets forms a $\sigma$-ideal. We prove an intermediate lemma first.
\begin{lem}\label{smallwitness} Let $A$ be Haar meager and $\e>0$. Then, there is a witness function whose range is
contained in an $\e$ ball
around the origin. \end{lem}
\begin{pf} Let $K$ be a compact set and $f: K \rightarrow X$ be a witness function for $A$. Let $U$ be an open subset of $X$
with diameter less than $\frac{\e}{2}$ such that $U \cap f(K) \neq \emptyset$. Let $y \in U \cap f(K)$ and $L =\overline {f^{-1}(U)}$.
Define $g:L \rightarrow X$ by letting $g(x)=f(x)-y$ for all
$x \in L$. This $g$ satisfies the required property. \end{pf}

\begin{thm}\label{sigmaideal} $\hm$ is a $\sigma$-ideal. \end{thm}
\begin{pf} That $\hm$ is hereditary follows from the definition of Haar meager. We will show that $\hm$ is closed under
countable union. To this end, let $A_1, A_2, \dots$ be Borel Haar meager sets.  For each $i$, let $K_i$ be a compact metric space and $g_i:K_i \rightarrow X$ be witness function for $A_i$.  By Lemma~\ref{smallwitness}, we 
may assume that $g_i(K_i)$ is contained in the ball centered at the origin with radius $2^{-i}$. Simply let $K =\prod_{i=1}^{\infty}K_i$. Now we define $g:K \rightarrow X$ as follows.
$g(x)=\sum_{i=1}^{\infty} g_i(x_i)$ where $x=(x_1,x_2,\ldots)$.  Clearly, $g$ is well-defined and continuous. To
complete the proof of the theorem, it will suffice to show
that $g$ is a witness for each $A_i$. To this end, fix an $i \in \N$ and $t \in X$. Note that $g^{-1}(A_i+t) = \{(x_1,
x_2, x_3,\ldots)\in \Pi_{i=1}^{\infty} K_i: \sum_{i=1}^{\infty} g_i(x_i) \in (A_i +t)\}$. Fix $x_1,x_2,\ldots,x_{i-1},x_{i+1}, \ldots$.
Since $g_i$ is a witness for $A_i$, we have that \[ \{a \in K_i: g_i(a) \in A_i +t - \sum_{k \neq i} g_k(x_k)
 \}\] is meager in $K_i$. As $A_i$ is Borel and $g$ is continuous, we have that $g^{-1}(A_i+t)$
 has the Baire property.  Applying the Kuratowski-Ulam Theorem, we obtain
that\\ $g^{-1}(A_i+t)$ is meager in $\Pi_{i=1}^{\infty} K_i$, completing the proof.
\end{pf}
We next show that in the definition of Haar meager set, one can choose $K \subseteq X$. However,
we prove a lemma first.
\begin{lem}\label{cantoroncompact} Let $K$ be a Cantor set and $M$ be any compact metric space. Then, there is
a continuous function $f$ from $K$ onto $M$ such that if $A$ is  meager in $M$ then $f^{-1}(A)$ is meager in
$K$. \end{lem}
\begin{pf} Using a standard construction from general topology, one can obtain a continuous function $f$ from $K$ onto $M$
such that if $U$ is nonempty and open in $K$ then $f(U)$ contains a nonempty open subset of $M$. We claim this $f$ has the desired property.
To obtain a contradiction, assume that there is a meager set $A \ss M$ such that $f^{-1}(A)$ is not meager. Let $F_1, F_2,\ldots$
be a sequence of nowhere dense closed sets such that $A \ss \cup_{i=1}^{\infty} F_i$. Since $f^{-1}(A)$ is not meager,
$f^{-1}(F_i)$ contains a nonempty open set for some $i$. This contradicts that the image of an open set under $f$ must contain
a nonempty open set.
\end{pf}
\begin{thm}\label{equivalentdef} Let $A \subseteq X$. Then the following are equivalent.
\begin{itemize}
\item Set $A$ is Haar meager.
\item There is a compact set $K \subseteq X$ and continuous function $g:K \rightarrow X$ which is a witness function for $A$.
\end{itemize}
\end{thm}
\begin{pf} Only ($\Rightarrow$) needs a proof as ($\Leftarrow$) is obvious. Let $A \subseteq X$ be a Borel Haar meager set. If $X$ is countable, then the only Haar meager subset of $X$ is the empty set and any set $K \subseteq X$ and $g:K \rightarrow X$ will do. If $X$ is uncountable, then it contains a copy of the Cantor space. Let $K$ be one such copy. Let $M$ be a compact metric space and $g:M \rightarrow X$ witness the fact that $A$ is Haar meager. Let $f:K \rightarrow M$ be a function
of the type in Lemma~\ref{cantoroncompact}. Then $g\circ f$ is the desired function. \end{pf}

\section{Remarks and Problems}
In this section we make some remarks, state some open problems concerning Haar meager sets and suggest possible lines of investigation which may be fruitful.

\begin{prob} Is the Continuum Hypothesis necessary in Example~\ref{ex1}? Can such an example be constructed in ZFC?
\end{prob}
Next problems concern the witness function in the definition of Haar meager set. 

\begin{prob}\label{one} Suppose $A \subseteq X$ is Haar meager. Is there a compact set $K \subseteq X$ such
that $(A+t) \cap K $ is meager in $K$ for all $t \in X$? 
\end{prob}

A negative answer to Problem~\ref{one}, would beg the following question.
\begin{prob}\label{two} Does the collection of all sets from Problem~\ref{one} form a $\sigma$-ideal?
\end{prob}

The next set of problems concerns concrete examples of Haar meager sets. As noted earlier, Kolar  \cite{kolar} defined the notion of HP-small set. This notion is finer than the notion of Haar null set and the notion of Haar meager set. As a corollary to his main result we obtain that the set of functions in $C([0,1])$ which have finite derivative at some point is HP-small. What we would like to do is distinguish the notion of Haar null set from the notion of Haar meager set.

\begin{prob}Are there concrete examples in analysis, topology, etc of sets which are Haar null but not Haar meager? Or vice versa? 
\end{prob}

As was observed earlier, every Polish group which is not locally compact admits  meager sets which are not Haar meager. Concrete examples of such sets in analysis are given by Zaj\'i\v cek
\cite{z} and Darji and White \cite{ds}. Zaj\'i\v cek \cite{z} showed that the set of functions in $C([0,1])$
for which $f'(0) = \infty$ has the property that it contains a translate of every compact subset of $C([0,1])$. Inspired by the result of Zaj\'i\v cek, Darji and White \cite{ds} gave an example of an uncountable, pairwise disjoint collection of subsets of $C([0,1])$ with the property that each element is meager yet contains a translate of every compact subset of $C([0,1])$.

\begin{prob} What are some natural examples in topology, dynamical systems, analysis of sets of natural interest which are meager but not Haar meager?
\end{prob}

We note that the notions of Haar meager is the topological mirror of the notion of Haar null. In Christensen's definition of Haar null, one can choose the test measure $\mu$ as a push forward of some
measure defined on some Cantor set in $X$. However, there are many instances in which these
two notions coincide.

\begin{prob} Are there some general conditions that one can put on a set so that Haar null $+$ the conditions implies Haar meager? The same question with the role of Haar meager and Haar null is reversed.
\end{prob}
Of course, Kolar's notion of HP-small implies both Haar null and Haar meager. What we are looking for some external condition of sets which allows transference principle. 
\section{acknowledgement} 
The author would like to thank the anonymous referee for valuable suggestions that improved the exposition of the paper.

\end{document}